\documentclass[11pt]{article}
\usepackage{amsthm}
\usepackage{amsmath,amsfonts,amssymb,rotating,colordvi}
\usepackage{float}

\usepackage[unicode=true,
 bookmarks=true,bookmarksnumbered=true,bookmarksopen=true,bookmarksopenlevel=2,
 breaklinks=false,pdfborder={0 0 1},backref=false,colorlinks=true]
 {hyperref}

\usepackage{tikz}

\usetikzlibrary{calc}
\usetikzlibrary{shapes.geometric}

\newcommand{\red}{\mathrm{red}}

\tikzset{
  c/.style={every coordinate/.try}
}

\usetikzlibrary{arrows,shapes,positioning}
\usetikzlibrary{decorations.markings}
\tikzstyle arrowstyle=[scale=1]
\tikzstyle directed=[postaction={decorate,decoration={markings,mark=at position 0.6 with {\arrow[arrowstyle]{stealth};}}}]
\tikzstyle reverse directed=[postaction={decorate,decoration={markings,mark=at position 0.4 with {\arrowreversed[arrowstyle]{stealth};}}}]
\tikzstyle dot=[style={circle,inner sep=1pt,fill}]
\hypersetup{pdftitle={Word-representability of triangulations of grid-covered cylinder graphs},
 pdfauthor={Zongqing Chen, Sergey Kitaev and Brian Y. Sun,...},
 pdfsubject={ },
 pdfkeywords={ },
 linkcolor=black,citecolor=black,urlcolor=blue,filecolor=blue,pdfpagelayout=OneColumn,pdfnewwindow=true,pdfstartview=XYZ,plainpages=false}

\def\qed{\nopagebreak\hfill{\rule{4pt}{7pt}}}

 \newtheorem{thm}{Theorem}[section]

\newtheorem{lem}[thm]{Lemma}

\newtheorem{coro}[thm]{Corollary}

\newtheorem{exa}[thm]{Example}

\numberwithin{equation}{section}

\newdimen\Squaresize \Squaresize=11pt
\newdimen\Thickness \Thickness=0.7pt
\def\Square#1{\hbox{\vrule width \Thickness
   \vbox to \Squaresize{\hrule height \Thickness\vss
    \hbox to \Squaresize{\hss#1\hss}
   \vss\hrule height\Thickness}
\unskip\vrule width \Thickness} \kern-\Thickness}

\def\Vsquare#1{\vbox{\Square{$#1$}}\kern-\Thickness}

\def\moins{\raise 1pt\hbox{{$\scriptstyle -$}}}

\begin{document}

\begin{center}
{\large \bf  On 132-representable graphs}
\end{center}

\begin{center}
Alice L.L. Gao$^{1}$,
Sergey Kitaev$^{2}$, and Philip B. Zhang$^{3}$\\[6pt]

$^{1}$Center for Combinatorics, LPMC-TJKLC \\
Nankai University, Tianjin 300071, P. R. China\\[6pt]

$^{2}$ Department of Computer and Information Sciences \\
University of Strathclyde, 26 Richmond Street, Glasgow G1 1XH, UK\\[6pt]

$^{3}$ College of Mathematical Sciences \\
Tianjin Normal University, Tianjin  300387, P. R. China\\[6pt]

Email: $^{1}${\tt gaolulublue@mail.nankai.edu.cn},
	   $^{2}${\tt sergey.kitaev@cis.strath.ac.uk},
           $^{3}${\tt zhangbiaonk@163.com}
\end{center}

\par
\noindent \textbf{Abstract.}
A graph $G = (V,E)$ is word-representable if there exists a word $w$ over the alphabet $V$ such that letters $x$ and $y$ alternate in $w$ if and only if $xy$ is an edge in $E$. Word-representable graphs are the subject of a long research line in the literature initiated in \cite{KP}, and they are the main focus in the recently published book \cite{KL}.  A word $w=w_1\cdots w_{n}$ avoids the pattern  $132$ if there are no $1\leq i_1<i_2<i_3\leq n$ such that $w_{i_1}<w_{i_3}<w_{i_2}$. The theory of patterns in words and permutations is a fast growing area discussed in \cite{HM,Kit}.

A research direction suggested in \cite{KL} is in merging the theories of word-representable graphs and patterns in words.  Namely, given a class of pattern-avoiding words, can we describe the class of graphs represented by the words? Our paper provides the first non-trivial results in this direction. We say that a graph is 132-representable if it can be represented by a 132-avoiding word. We show that each 132-representable graph is necessarily a circle graph. Also, we show that any tree and any cycle graph are 132-representable, which is a rather surprising fact taking into account that most of these graphs are non-representable in the sense specified, as a generalization of the notion of a word-representable graph, in \cite{JKPR}. Finally, we provide explicit 132-avoiding representations for all graphs on at most five vertices, and also describe all such representations, and enumerate them, for complete graphs.\\


\section{Introduction}
A graph $G = (V,E)$ is word-representable if there exists a word $w$ over the alphabet $V$ such that letters $x$ and $y$ alternate in $w$ if and only if $xy$ is an edge in $E$. For example, the graph to the right in Figure~\ref{circle-graph} is word-representable and one of words representing it is $bcdad$. Some graphs are word-representable, others are not, and the minimum non-word-representable graph is the wheel $W_5$ shown to the left in Figure~\ref{wheel5}.

Word-representable graphs are the subject of a long line of research in the literature initiated in~\cite{KP}, and they are the main focus in the recently published book \cite{KL}. A general program of research suggested in \cite[p.~183]{KL} takes as the input a language defined, for example, through pattern-avoiding words, and outputs a description of the class of graphs represented by the language. For instance, as is discussed in \cite[p.~183]{KL}, the set of weakly increasing words (those {\em avoiding the pattern} 21) defines graphs whose vertices can be partitioned into a clique and an independent set, so that no edge connects the clique and the independent set. However, apart from this simple result, no research has been done in this direction.

In this paper, we study graphs defined by 132-avoiding words. Our research merges the theories of word-representable graphs  \cite{KL} and patterns in words \cite{HM,Kit}, the latter being a very fast growing area. A word $w=w_1w_2\cdots w_{n}$
avoids the pattern $132$ (resp., $123$) if there are no indices $1\leq i_1<i_2<i_3\leq n$ such that $w_{i_1}<w_{i_3}<w_{i_2}$ (resp., $w_{i_1}<w_{i_2}<w_{i_3}$). We say that a graph $G$ is 132-representable (resp., 123-representable) if there is a 132-avoiding (resp., 123-avoiding) word representing it. Note that for the last definition to make sense, labels of graphs are supposed to be taken from a totally ordered set. Also, when trying to 132-represent (123-represent) a graph, we are allowed to label the graph in any suitable way\footnote{There is no issue with labelling when considering word-representable graphs, since all labelings are equally good or bad. However, in the contexts when there is an order on labels, labelling graphs in a proper way may be essential for finding a representation.}.

\begin{figure}[H]
    \centering
    \begin{minipage}{.95\textwidth}
        \centering
        \begin{tikzpicture}[scale=1]

\draw[thin] (0,0) ellipse (65pt and 20pt);
\draw[thin] (0,0.3) ellipse (75pt and 35pt);
\draw[thin] (0,0.7) ellipse (100pt and 50pt);
\draw[thin] (0,1.0) ellipse (110pt and 60pt);

\node  at (0,0) { $132$-representable~graphs };
 \node  at (0,1.0) {circle~graphs};
 \node  at (0,1.9) {word-representable~graphs};
 \node at (0,2.8) { all~graphs };

\draw [thick, <-] (1.3,2.8) -- (3.3,2.8); \node[right] at (3.3,2.8) {odd wheels $W_5$, $W_7,\ldots$ \cite{KP}};
\draw [thick, <-] (1.5,1.5) -- (4,1.5); \node[right] at (4,1.5) {prisms \cite{Kit1}};
\draw [thick, <-] (1.5,0.8) -- (4,0.8); \node[right] at (4,0.8) {disjoint union of 2 comp-};
															\node[right] at (4,0.4) {lete graphs of size 4 \cite{M}};
\draw [thick, <-] (1,-0.4) -- (3.3,-0.4); \node[right] at (3.3,-0.4) {trees, cycle graphs,};  \node[right] at (3.3,-0.8) {\ \ complete graphs};

  \end{tikzpicture}
    \end{minipage}
  \caption{The place of $132$-representable graphs in a hierarchy of graph classes}
  \label{graph=class}
\end{figure}
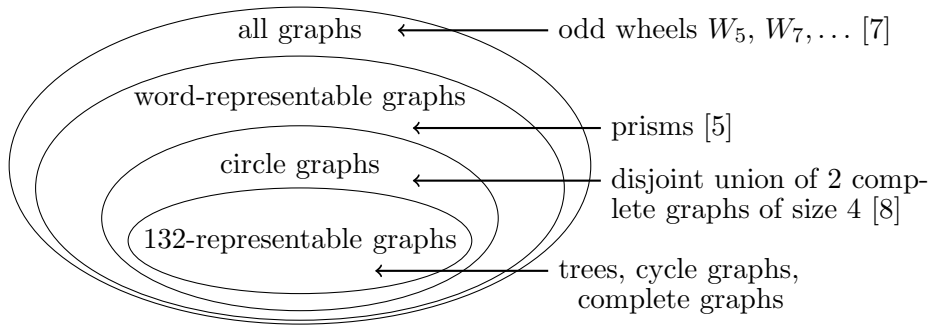

One of the main results in this paper is in showing that any 132-represent-able graph is necessarily a circle graph. A result in \cite{M} shows that 132-representable graphs are a strict subset of circle graphs. Also, we show that trees, cycle graphs and complete graphs are 132-representable. Thus, the place of 132-representable graphs in a hierarchy of graph classes is as shown in Figure~\ref{graph=class}, where we also indicate known facts  that odd wheels are non-word-representable \cite{KP}, while prisms are word-representable but not circle graphs \cite{Kit1}. Interestingly, the studies in \cite{M} show that the class of 123-representable graphs, being different from the class of 132-representable graphs, is also a proper subclass of circle graphs, even though not all trees are 123-representable; all cycle graphs and complete graphs are 123-representable.

One should compare our results with the results on 12-{\em representable graphs} obtained in~\cite{JKPR}. These graphs are an instance of $u$-representable graphs, a far reaching generalization of word-representable graphs, also introduced in~\cite{JKPR}, where $u$ is a word over $\{1,2\}$ different from $22\cdots 2$. Similarly to the case of 132-representable graphs, labelling of graphs is important for 12-representation. A word $w$ 12-represents a graph $G$, if for any labels $x$ and $y$, $x<y$, $xy$ is an edge in $G$ if and only if after removing all letters in $w$ but $x$ and $y$, we will obtain a word of the form $yy\cdots yxx\cdots x$. Note that the notions of 132-representable graphs and 12-representable graphs are not directly related (in the former case the pattern is used to give a condition on words representing graphs, while in the latter case the pattern is used to define the representation itself). It was shown in \cite{JKPR} that any 12-representable graph is necessarily a comparability graph, while very few trees (called {\em double caterpillars}) and almost no cycle graphs (only cycle graphs on at most four vertices) are 12-representable.

This paper is organized as follows. In Section~\ref{sec2} we give necessary definitions, notation and results to be used in the paper. In Section~\ref{132-representants-sec} we  derive a key property of words 132-representing graphs (see Theorem~\ref{thm3-2}) and state its corollary, the main result in this paper, that any 132-representable graph is necessarily a circle graph (see Corollary~\ref{thm-circle}). In Section~\ref{sec-misc} we not only establish 132-representability of trees and cycle graphs, but also describe and enumerate all 132-representants for complete graphs. Moreover, in Section~\ref{sec-misc} we discuss non-132-representable graphs and give explicit 132-representation of graphs on four and five vertices. Finally, in Section~\ref{open-sec} we state a number of suggestions for further research.

\section{Preliminaries}\label{sec2}

\noindent
{\bf Graphs.} We will now review a number of basic notions/notations in graph theory. In this paper, we deal with  \emph{simple graphs}, that is, graphs with no {\em loops} and no {\em multiple edges}.

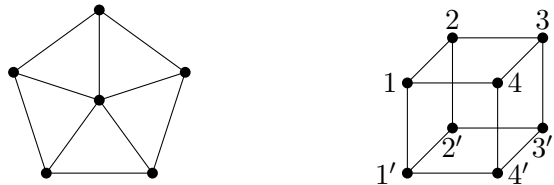
\begin{figure}[H]
    \centering
    \begin{minipage}{.95\textwidth}
        \centering
        \begin{tikzpicture}[scale=1.2]

       \draw (0,0)--(-0.363,1.118)--(0.587,1.809)--(1.539,1.118)--(1.176,0)--(0,0);
       \draw (0.587,0.809)--(0,0);
       \draw (0.587,0.809)--(-0.363,1.118);
       \draw (0.587,0.809)--(0.587,1.809);
       \draw (0.587,0.809)--(1.539,1.118);
       \draw (0.587,0.809)--(1.176,0);

       \fill(0,0) circle(0.06cm);\fill(1.176,0) circle(0.06cm);
       \fill(-0.363,1.118) circle(0.06cm);\fill(1.539,1.118) circle(0.06cm);
       \fill(0.587,0.809) circle(0.06cm); \fill(0.587,1.809) circle(0.06cm);

 \draw (4,0)--(5,0)--(5,1)--(4,1)--(4,0);
       \draw (4.5,0.5)--(5.5,0.5)--(5.5,1.5)--(4.5,1.5)--(4.5,0.5);
       \draw (4,0)--(4.5,0.5);
       \draw (5,0)--(5.5,0.5);
       \draw (5,1)--(5.5,1.5);
       \draw (4,1)--(4.5,1.5);

       \fill(4,0) circle(0.06cm);\fill(5,0) circle(0.06cm);
       \fill(5,1) circle(0.06cm);\fill(4,1) circle(0.06cm);
       \fill(4.5,0.5) circle(0.06cm); \fill(5.5,0.5) circle(0.06cm);
       \fill(5.5,1.5) circle(0.06cm);\fill(4.5,1.5) circle(0.06cm);

       \node [left] at (4,0) {  $1'$}; \node [left] at (4,1) {  $1$};
       \node  at (4.5,0.3) {  $2'$};\node at (4.5,1.7) {  $2$};
       \node   at (5.5,0.3) {  $3'$};\node at (5.5,1.7) {  $3$};
       \node [right] at (5,0) {  $4'$};\node [right] at (5,1) {  $4$};

        \end{tikzpicture}

    \end{minipage}
  \caption{The wheel graph $W_5$ and the prism $Pr_4$}
  \label{wheel5}
\end{figure}

The \emph{degree} $d(v)$ of a vertex $v$ in a graph $G$ is the number
of edges of $G$ incident with $v$. The {\em complete graph} on $n$ vertices is denoted by $K_n$.
A {\em cycle graph} $C_n$ is the graph on $n$ vertices that consists of a single cycle.  A {\em wheel graph} $W_n$ is the graph on $n+1$ vertices obtained from
$C_n$ by adding an all-adjacent vertex ({\em apex}). The wheel graph $W_5$ is shown to the left in Figure~\ref{wheel5}.

%
%
%
%

A {\em prism} $Pr_n$ is a graph consisting of two cycles $12\cdots n$ and $1'2'\cdots n'$, where $n\geq 3$, connected by the edges $ii'$ for $i=1,\ldots,n$. For example, $Pr_4$, also known as the {\em three-dimensional cube}, is shown to the right in Figure~\ref{wheel5}.

\begin{figure}[H]
    \centering
    \begin{minipage}{.95\textwidth}
        \centering
        \begin{tikzpicture}[scale=1]

\draw[thin] (0,0) circle (1cm);
\draw (-1,0)--(1,0);
\draw (0,1)--(0,-1);
\draw (-0.707,-0.707)--(0.707,0.707);
\draw (-0.6,-0.8)--(0.6,-0.8);

\fill(-1,0) circle(0.06cm);
\fill(1,0) circle(0.06cm);
\fill(0,1) circle(0.06cm);
\fill(0,-1) circle(0.06cm);
\fill(-0.6,-0.8) circle(0.06cm);
\fill(0.6,-0.8) circle(0.06cm);
\fill(-0.707,-0.707) circle(0.06cm);
\fill(0.707,0.707) circle(0.06cm);

 \node  at (-0.2,0.6) {$a$};
 \node [right] at (0.1,0.65) {$b$};
 \node [right] at (0.4,0.2) {$c$};
 \node [right] at (0.1,-0.55) {$d$};

 \node  at (3,0) {  $\Leftrightarrow$};

\draw (5,-0.5)--(6,-0.5)--(7,-0.5);
\draw (5,-0.5)--(6,0.5)--(6,-0.5);

\fill(5,-0.5) circle(0.06cm); \fill(6,-0.5) circle(0.06cm);
\fill(7,-0.5) circle(0.06cm); \fill(6,0.5) circle(0.06cm);

\node at (5,-0.7) {$c$};
\node at (6,-0.7) {$a$};
\node at (7.1,-0.7) {$d$};
\node at (6,0.8) {$b$};

        \end{tikzpicture}
    \end{minipage}
  \caption{A circle with four chords and the corresponding circle graph}
  \label{circle-graph}
\end{figure}
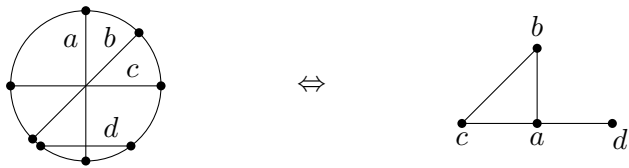

Finally, a \emph{circle graph} is a graph whose vertices can be associated with chords of a circle such that two vertices are adjacent if and only if the corresponding chords intersect. See Figure~\ref{circle-graph} for an example of a circle graph.

\noindent
{\bf Words and permutations.} For a finite word $w$, let $A(w)$ denote the set of letters occurring
in~$w$, and $\red(w)$ denote the word over $\{1,2,\ldots,|A(w)|\}$ obtained by replacing the $i$-th smallest letter(s) by $i$. We call $\red(w)$ the {\em reduced form} of $w$. Also, for any $x\in A(w)$, let $n_w(x)$ denote the number of copies of $x$ in $w$, and $x_i$ denote the $i$-th occurrence of $x$ in $w$ from left to right. For example, if $w=14661476212$, then $A(w)=\{1,2,4,6,7\}$, $\red(w)=13441354212$, and say for $x=6$, $n_w(6)=3$. A word $w$ is {\em $k$-uniform} if each letter in $w$ occurs exactly $k$ times.

Suppose that $x$ and $y$ are two distinct letters in $A(w)$.
We say that $x$ and $y$ {\it alternate} in $w$ if after deleting in $w$ all letters but the copies of $x$ and $y$ we either obtain a word
$xyxy\cdots$ (of even or odd length) or a word $yxyx\cdots$ (of even or odd length).
In particular, if $w$ has a single occurrence of $x$ and a single occurrence of $y$, then $x$ and $y$ alternate in $w$.

A word or permutation $w=w_1w_2\cdots w_{n}$ {\em avoids} the pattern $132$ if there are no indices $1\leq i_1<i_2<i_3\leq n$ such that $w_{i_1}<w_{i_3}<w_{i_2}$. For example, the word $31458$ avoids the pattern 132, while $3474$ is not 132-avoiding (the subsequence 374 in this word forms the pattern 132).  It is  a well-known fact (e.g. see~\cite[p.~32]{Kit}) that the number of 132-avoiding permutations of length $n$ is given by the {\em $n$-th Catalan number} $C_n=\frac{1}{n+1}{2n\choose n}$.

A subword of $w$ formed by consecutive letters is called a {\em factor} of $w$. For example, 6651 and 41 are factors of $26651141$.  Finally, we let $[n]=\{1,2,\ldots,n\}$.\\

\noindent
{\bf Word-representable graphs.} A graph $G=(V,E)$ is \emph{word-representable} if there exists a word $w$ over the alphabet $A(w)=V$ such that $x$ and $y$ alternate in $w$ if and only if $xy\in E$ for each $x\neq y$ (that is, $x$ and $y$ are connected by an edge).
In this context, we say that $w$ {\em represents} $G$ and $w$ is a {\em word-representant} for $G$.

In this paper we assume that elements in $V$ come from a totally ordered alphabet, which is important for the following definition. A word-representable graph $G$ is {\em $132$-representable} if, possibly after relabelling the graph, there exists a $132$-avoiding word $w$ that represents $G$. In this context, $w$ is called a {\em $132$-representant} for $G$.

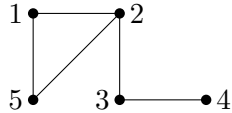
\begin{figure}[!htb]
    \centering
    \begin{minipage}{.95\textwidth}
        \centering
        \begin{tikzpicture}[scale=1.15]

     \draw (0,0)--(1,1);
     \draw (0,0)--(0,1)--(1,1)--(1,0)--(2,0);

     \fill(0,0) circle(0.06cm);\fill(1,0) circle(0.06cm);
     \fill(1,1) circle(0.06cm);\fill(0,1) circle(0.06cm);
     \fill(2,0) circle(0.06cm);

       \node  [left] at (0,0) {$5$};\node  [right] at (0.6,0) {$3$};
       \node  [left] at (0,1) {$1$};\node   [right] at (1,1) {$2$};
       \node [right] at (2,0) {$4$};

\end{tikzpicture}
\end{minipage}
 \caption{A $132$-representable graph $G$}
  \label{wordandgraph}
\end{figure}

For example, if $w=43451251$, then the subword induced by the letters $1$ and $2$
is $121$, and hence the letters $1$ and $2$ alternate in $w$, so that the respective vertices are connected in $G$. On the other hand, the letters $1$ and $3$
do not alternate in $w$, because removing all other letters we obtain $311$; thus, $1$ and $3$ are not connected in $G$. Figure~\ref{wordandgraph} shows the graph represented by $w$.
Moreover, since $w$ is $132$-avoiding, $G$ is $132$-representable and $w$ is a $132$-representant of $G$.

\begin{figure}[!htb]
    \centering
    \begin{minipage}{.95\textwidth}
        \centering
        \begin{tikzpicture}[scale=1.15]

    \draw (-3,0)--(-2,0)--(-1,0); \draw (-2,0)--(-2,1);
     \fill(-3,0) circle(0.06cm);\fill(-2,0) circle(0.06cm);  \fill(-1,0) circle(0.06cm);\fill(-2,1) circle(0.06cm);

    \node [left] at (-1.6,-0.4) {$A$.};

    \draw (0,0)--(1,0)--(2,0); \draw (1,0)--(1,1);
     \fill(0,0) circle(0.06cm);\fill(1,0) circle(0.06cm);  \fill(2,0) circle(0.06cm);\fill(1,1) circle(0.06cm);
      \node  [left] at (0,0) {$2$};\node  [right] at (1,0.2) {$1$};
       \node  [left] at (2.4,0) {$4$};\node   [right] at (1,1) {$3$};

   \node [left] at (1.4,-0.4) {$B$.};

    \draw (3,0)--(4,0)--(5,0); \draw (4,0)--(4,1);
     \fill(3,0) circle(0.06cm);\fill(4,0) circle(0.06cm);  \fill(5,0) circle(0.06cm);\fill(4,1) circle(0.06cm);
         \node  [left] at (3,0) {$2$};\node  [right] at (4,0.2) {$4$};
       \node  [left] at (5.4,0) {$1$};\node   [right] at (4,1) {$3$};

 \node [left] at (4.4,-0.4) {$C$.};

\end{tikzpicture}
\end{minipage}
 \caption{Significance of proper labelling}
  \label{proper-label-ex}
\end{figure}
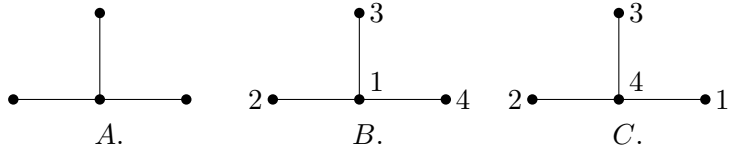

We note that labelling of a graph is important when dealing with $132$-representation, which is not the case with just word-representation since all labellings are equally good or bad. For example, the fact that the (unlabelled) graph $A$ in Figure~\ref{proper-label-ex} is 132-representable is given by the labelled version $B$ of it and the 132-avoiding word 43212341. However, if we would label $A$ to obtain the graph $C$  in Figure~\ref{proper-label-ex}, then no 132-avoiding representation of it exists. Indeed, suppose that a 132-representant $w$ for $C$ exists. Then at least two letters in $\{1,2,3\}$, say $x$ and $y$, $x<y$, must be repeated at least twice in $w$, or else there would be at least one unwanted edge in $\{12, 13, 23\}$. Further, because 4 is an apex, there are $x$'s and $y$'s on both sides of a 4 in $w$ (the 4 must alternate with $x$ and $y$), which leads to an occurrence $x4y$ of the pattern 132; contradiction.

The following result is of special importance to us.

\begin{thm}[\cite{Halldorsson15}]\label{circle-thm} A graph $G$ is word-representable and its representation requires at most
two copies of each letter if and only if $G$ is a circle graph. \end{thm}

Another relevant result is as follows.

\begin{thm}[\cite{Kit1}]\label{3-repr-not-2-repr} Prisms cannot be represented using at most two copies of each letter, but can be represented using at most three copies of each letter. \end{thm}

\section{$132$-representants}\label{132-representants-sec}

In this section, we discuss some properties of $132$-representants.

We first present a simple, but useful theorem.

\begin{thm}\label{degree2}
Let $G$ be a $132$-representable graph, and $x$ be a vertex in $G$ such that $d(x)\geq 2$. Then for any $132$-representant $w$ of $G$, we have $n_w(x)\leq 2$.
\end{thm}

\proof
Since $d(x)\geq 2$, there exist vertices $a$ and $b$, $a>b$, in $G$ that are adjacent with $x$.

Suppose that there are at least three copies of $x$ in $w$. Then by the definition of a $132$-representant,
there exists a subsequence $xw_1xw_2x$ in $w$, where for $i=1,2$, $w_i$ is a factor of $w$ containing exactly one $a$,  one $b$, and no $x$.
There are three cases to consider, all of which contradict the requirement that $w$ is $132$-avoiding:
\begin{itemize}
\item $x>a>b$:  $bxa$ is a 132 pattern in $w$ where $b\in w_1$ and $a\in w_2$;
\item $a>b>x$:  $xab$ is a 132 pattern in $w$ where $a\in w_1$ and $b\in w_2$;
\item $a>x>b$:  $bax$ is a 132 pattern in $w$ where $b\in w_1$ and $a\in w_2$.
\end{itemize}
Hence, at most two copies of $x$ can appear in $w$.
\qed

As consequences of Theorem~\ref{degree2}, we obtain the following results.
\begin{coro}
If each vertex in a graph $G$ is of degree at least $2$, then any $132$-representant for $G$ is of length at most $2n$.
\end{coro}

\begin{coro}\label{degree1-3}
Let $w$ be a $132$-representant for a graph $G$. If $d(x)=1$ and the vertex $a$ connected to $x$ has degree at least $2$, then $x$ occurs at most three times in $w$.
\end{coro}

\proof
Let $w$ denote a 132-representant for $G$.
Since $d(a)\geq 2$, by Theorem~\ref{degree2} $a$ occurs at most twice in $w$.
Combining with the fact that $a$ and $x$ alternate in $w$, we have that $x$ occurs at most three times in $w$.
\qed

The following theorem generalizes Theorem~\ref{degree2}.

\begin{thm}\label{thm3-2}
If a graph $G$ is $132$-representable, then there exists a $132$-avoiding word $w$ representing $G$ such that for any letter $x$ in $w$,
$n_w(x)\leq 2.$
\end{thm}

\proof
Let $w$ be a 132-representant for $G$.
If all the vertices in $G$ have degree at least 2, then by Theorem \ref{degree2} every letter appears in $w$ at most twice.
Hence it suffices to consider the case where there exists a vertex $x$ in $G$ such that $d(x)=1$. Let $a$ be the vertex connected to $x$. We consider two cases.

\begin{itemize}
\item $d(a)\geq2$.
By Corollary \ref{degree1-3}, the letter $x$ occurs at most three times in $w$.
To prove the theorem, we assume that there are three copies of $x$ in $w$ and then we will  construct a new $132$-avoiding word $w'$ which also represents $G$ but contains only two copies of  $x$. By
Theorem \ref{degree2}, there are exactly two copies of $a$ in $w$. In what follows, according to our notation, $x_i$ denotes the $i$-th $x$ and $a_j$ the $j$-th $a$ in $w$ from left to right, where $1\leq i\leq 3$ and $1\leq j\leq 2$.

Suppose that $a>x$. If there are no letters between the $a$'s except for $x$ then $a$ is connected only to $x$ in $G$; contradiction with $d(a)\geq 2$. Thus there is a letter $b\neq x$ between $a_1$ and $a_2$ in $w$.  If $b>a>x$, then $x_1ba_2$ will be the pattern 132; if $a>b>x$, then $x_1a_1b$ will form the pattern 132; if $a>x>b$, then $ba_2x_3$ will form the pattern 132; in either case, there is a contradiction with the definition of $w$. Thus we must have $a<x$.

We next construct a new $132$-avoiding word $w'$ from $w$.
Since there is no element $t$ smaller than $a$ to the left of $a_1$ in $w$ (or else, $tx_2a_2$ would be the $132$-pattern), we obtain that $a$ is a left-to-right minimum in $w$ (that is, no letter to the left of $a$ is less than $a$). We delete all three $x$'s and replace $a_1$ by the factor $a^+a_1a^+$ to obtain the new word $w'$, where $a<a^+<a+1$. By construction of $w'$, if it contains an occurrence of the pattern $132$ then this occurrence cannot involve $a^+$ and thus it would give an occurrence of the pattern in $w$; contradiction. Moreover, $a$ is the only letter in $w'$ alternating with $a^+$, and thus $w'$ 132-represents $G'$ obtained from $G$ by replacing the label  $x$ by $a^+$.

\item $d(a)=1$, which means that the edge $xa$ is disconnected from the rest of the graph. Let $w'$ denote the word obtained from $w$ by deleting $a$ and $x$. Clearly, $w'$ is 132-avoiding.  But then the 132-avoiding word $n(n-1)n(n-1)w'$, where $n$ and $n-1$ are larger than any other letter in $A(w')$, represents the graph $G'$ obtained from $G$ by replacing the labels $a$ and $x$ by $n$ and $n-1$ (in any order).
\end{itemize}
We can repeat the procedure described above for any other vertices of degree~1 in $G$ to obtain the desired result.
\qed

One of the main results in this paper is the following statement.

\begin{coro}\label{thm-circle}
Any $132$-representable graph is a circle graph.
\end{coro}

\proof
Let $G$ be a $132$-representable graph. By Theorem \ref{thm3-2}, there exists a $132$-representant $w$ of $G$ that contains at most two copies of each letter. By Theorem~\ref{circle-thm} $G$ is a circle graph.\qed

Note that we do not know whether each circle graph is 132-representable or not.

\section{$132$-representable graphs}\label{sec-misc}

In this section, we will show that trees, cycles, and complete graphs are 132-representable.

\subsection{Trees and cycle graphs}

\begin{thm}\label{tree}
Trees are $132$-representable.
\end{thm}

\proof
We proceed by induction on the number of vertices with an additional condition. The tree with only one vertex can be represented by $1$. Suppose that we can represent a tree with less than $n$ vertices by a 132-avoiding word and the label of the root has only one occurrence and the label of the non-root vertex has exactly two occurs in the corresponding word.

Given a tree $T$ with $n$ vertices, label it in pre-order, that is, starting from the root traverse the subtrees from left to right recursively. See the graph to the left in Figure~\ref{treeA} for an example. Suppose that the root has $r$ children, which means that $T$ has $r$ subtrees, whose roots are children of the root of $T$. Denote the $r$ trees by $T_i$ for $1\leq i \leq r$ from left to right and suppose that the root of $T_i$ is labeled by $n_i$. Note that $2\leq n_1< n_2<\cdots< n_r\leq n$, so that for $1\leq i\leq r$, $T_i$ has $n_{i+1}-n_i$ vertices, where $n_{r+1}=n+1$. Hence $T_i$ is a tree having less than $n$ vertices.
By induction hypothesis, $T_i$ is 132-representable and it can be represented by a 132-avoiding word $w(T_i)$ with only one copy of $n_i$ and two copies of any other letter. Let $w=w(T_r)w(T_{r-1})\cdots w(T_1)1n_1n_2\cdots n_r$. It is easy to see that $w$ represents $T$, and in particular, the root labeled by 1 is only connected to its children. Moreover, since for $1\leq i<j\leq r$ the labels of $T_i$ are smaller that the labels of $w(T_j)$, we get that $w$ is 132-avoiding. We are done.
\qed

\begin{exa}\textnormal{
Let $T$ be a tree as follows. It is clearly that $T$ has three subtrees $T_1$, $T_2$ and $T_3$. By Theorem \ref{tree}, there is $w(T_2)=5$. Moreover, we have
$w(T_1)=43234$ and $w(T_3)=87678$, which can be obtained by applying the inductive argument again. Hence $w(T)=87678.5.43234.1256$, where the dots showing parts of $w(T)$ should be ignored. It is obvious that $w(T)$ is $132$-avoiding and it represents $T$.}

\begin{figure}[!htb]
    \centering
    \begin{minipage}{.95\textwidth}
        \centering
        \begin{tikzpicture}[scale=1.5]

       \draw (2.1,0)--(2.5,0.5)--(2.9,0);
       \draw (2.5,0.5)--(3,1)--(3.5,0.5);
       \draw (3.1,0)--(3.5,0.5)--(3.9,0);
       \draw (3,0.5)--(3,1);

       \fill(2.1,0) circle(0.06cm); \fill(2.9,0) circle(0.06cm);
       \fill(3.1,0) circle(0.06cm); \fill(3.9,0) circle(0.06cm);
       \fill(2.5,0.5) circle(0.06cm);\fill(3,0.5) circle(0.06cm);
       \fill(3.5,0.5) circle(0.06cm);
       \fill(3,1) circle(0.06cm);

        \node [right] at (3,1) {  $1$};
       \node [right] at (2.5,0.5) {  $2$};
       \node [right] at (3,0.5) {  $5$};
       \node [right] at (3.5,0.5) {  $6$};
       \node  at (2.1,-0.2) {  $3$};
       \node  at (2.9,-0.2) {  $4$};
       \node  at (3.1,-0.2) {  $7$};
       \node  at (3.9,-0.2) {  $8$};
       \node  at (3,-0.5) {  $T$};

       \node  at (4.5,0.3) {  $\Rightarrow$};

       \draw (5.1,0)--(5.5,0.5)--(5.9,0);
       \fill(5.1,0) circle(0.06cm); \fill(5.9,0) circle(0.06cm);
       \fill(5.5,0.5) circle(0.06cm);
       \node [right] at (5.5,0.5) {  $2$};
        \node  at (5.1,-0.2) {  $3$};
       \node  at (5.9,-0.2) {  $4$};
       \node  at (5.5,-0.5) {  $T_1$};

       \fill(6.5,0.5) circle(0.06cm);
       \node [right] at (6.5,0.5) {  $5$};
       \node  at (6.5,-0.5) {  $T_2$};

       \draw (7.1,0)--(7.5,0.5)--(7.9,0);
       \fill(7.1,0) circle(0.06cm); \fill(7.9,0) circle(0.06cm);
       \fill(7.5,0.5) circle(0.06cm);
        \node [right] at (7.5,0.5) {  $6$};
        \node  at (7.1,-0.2) {  $7$};
       \node  at (7.9,-0.2) {  $8$};
      \node  at (7.5,-0.5) {  $T_3$};

        \end{tikzpicture}
    \end{minipage}
  \caption{A tree $T$ of size $8$ and its subtrees}
  \label{treeA}
\end{figure}
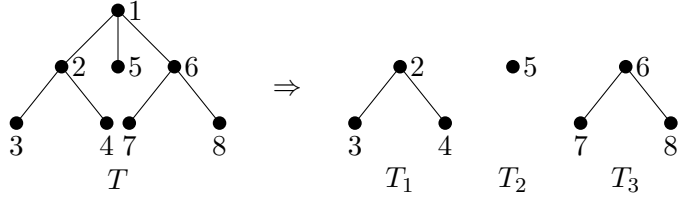
\end{exa}

\begin{figure}[!htb]
    \centering
    \begin{minipage}{.5\textwidth}
        \centering
        \begin{tikzpicture}[scale=1.5]

       \draw (5,0)--(6,0)--(7,0)--(7.2,0);
       \draw (7.8,0)--(8,0)--(9,0);

       \fill(5,0) circle(0.06cm); \fill(6,0) circle(0.06cm);
       \fill(7,0) circle(0.06cm); \fill(8,0) circle(0.06cm);
       \fill(9,0) circle(0.06cm);


\node [right] at (7.25,0) {  $\cdots$};

       \node [right] at (4.8,-0.2) {  $1$};
       \node [right] at (5.8,-0.2) {  $2$};
       \node [right] at (6.8,-0.2) {  $3$};
       \node [right] at (7.8,-0.2) {  $n-1$};
       \node [right] at (8.8,-0.2) {  $n$};
        \end{tikzpicture}
    \end{minipage}
  \caption{A path graph $P_n$}
  \label{path}
\end{figure}
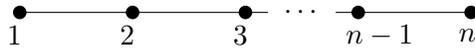

\begin{coro}\label{C-n-corol}
Cycle graphs $C_n$ are $132$-representable.
\end{coro}
\proof
Let $n\geq3$. A path graph $P_n$ (see Figure~\ref{path}) is a tree, and, by the proof of Theorem~\ref{tree}, it can be represented by the $132$-avoiding word
$$w=n(n-1)n(n-2)(n-1)(n-3)(n-2)\cdots 45342312.$$
Let $w'$ be the word obtained from $w$ by deleting the first $n$ in $w$. Then it is easy to see that  $w'$ represents $C_n$.
\qed

\begin{exa}\textnormal{
$132$-representants for $C_4$ and $C_5$, based on Corollary~\ref{C-n-corol}, are given in Figure \ref{path1}.}
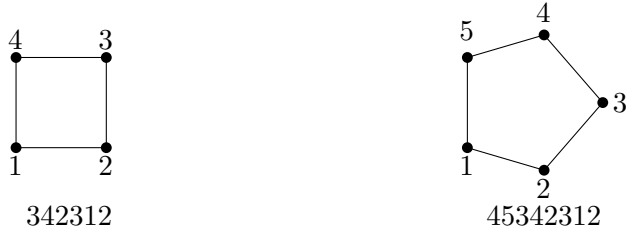
\begin{figure}[H]
    \centering
    \begin{minipage}{.95\textwidth}
        \centering
        \begin{tikzpicture}[scale=1.2]

       \draw (1,0)--(2,0)--(2,1)--(1,1)--(1,0);

 \fill(1,0) circle(0.06cm);\fill(2,0) circle(0.06cm);\fill(2,1) circle(0.06cm);\fill(1,1) circle(0.06cm);

       \node [right] at (0.8,-0.2) {  $1$};
       \node [right] at (1.8,-0.2) {  $2$};
       \node [right] at (1.8,1.2) {  $3$};
       \node [right] at (0.8,1.2) {  $4$};
       \node [right] at (1,-0.75) {  $342312$};

       \draw (6,0)--(6.85,-0.25)--(7.5,0.5)--(6.85,1.25)--(6,1)--(6,0);

 \fill(6,0) circle(0.06cm);\fill(6.85,-0.25) circle(0.06cm);\fill(7.5,0.5) circle(0.06cm);\fill(6.85,1.25) circle(0.06cm); \fill(6,1) circle(0.06cm);

       \node [right] at (5.8,-0.2) {  $1$};
       \node [right] at (6.65,-0.45) {  $2$};
       \node [right] at (7.5,0.5) {  $3$};
       \node [right] at (6.65,1.5) {  $4$};
       \node [right] at (5.8,1.3) {  $5$};
       \node [right] at (6.1,-0.75) {  $45342312$};
        \end{tikzpicture}
    \end{minipage}
  \caption{132-representants for $C_4$ and $C_5$}
  \label{path1}
\end{figure}

\end{exa}

\subsection{Complete graphs}
In the following theorem we shall describe and enumerate all 132-representants for $K_n$.

\begin{thm}\label{thm-K-n-132}
For $n\geq1$, a complete graph $K_n$ is $132$-representable. Moreover, for $n\geq 3$, there are $$2+C_{n-2}+\sum_{i=0}^{n}C_i$$ different $132$-representants for $K_n$, where $C_n=\frac{1}{n+1}{2n\choose n}$ is the $n$-th Catalan number. Finally, $K_1$ can be represented by a word of the form $11\cdots 1$ and $K_2$ by a word of the form $1212\cdots$ (of even  or odd length) or $2121\cdots$ (of even or odd length).
\end{thm}

\proof
Clearly, $K_1$ can only be represented by a word of the form $11\cdots 1$, and $K_2$ can only be represented by a word of the form $1212\cdots$ (of even  or odd length) or $2121\cdots$ (of even or odd length). Each of these words is 132-avoiding.

Let $n\geq3$. Suppose that $w$ is a 132-representant for $K_n$. According to the definition of a complete graph, for any $1\leq i<j\leq n$, we have that $i$ and $j$ alternate in $w$. Since $d(n)\geq 2$, by Theorem~\ref{degree2}, there are two cases to consider.\\

\noindent
{\bf Case 1.} There are exactly two copies of $n$ in $w$, and $w=w_1nw_2nw_3$, where $w_k$ is a word over $[n-1]$ for $k=1,2,3$.
    Since for $1\leq i\leq n-1$, $i$ and $n$ alternate in $w$, there is exactly 1 copy of $i$ in $w_2$, which means that $w_2$ is in fact a permutation of length $n-1$. Moreover, for $1\leq i\leq n-2$, $i$ must not appear in $w_1$, or $i,n,n-1$ will form the pattern 132. Thus, $w_1=n-1$ or $w_1=\epsilon$, the empty word. Similarly, we have that $w_3=1$ or $w_3=\epsilon$. Thus, there are four subcases to consider and in each subcase, we just need to consider the form of $w_2$.

\noindent
{\bf Subcase 1.1.} $w_1=n-1$ and $w_3=1$. Thus $1$ is to the left of $n-1$ in $w_2$, since $1$ and $n-1$ alternate in $w$. For $2\leq i\leq n-2$, $i$ must be between $1$ and $n-1$ in $w_2$ since $i$ alternates with $1$ and $n-1$. Moreover, for $2\leq i< j\leq n-2$, they are in increasing order in $w_2$, or $1,j,i$ will form a 132 pattern. Hence, we obtain that $w=(n-1)nw'1$ where $w'$ is the increasing permutation $12\cdots n$, and this case contributes one representation.

\noindent
{\bf Subcase 1.2.} $w_1=n-1$ and $w_3=\epsilon$. For $1\leq i\leq n-2$, $i$ is to the left of $n-1$ in $w_2$, since $i$ and $n-1$ alternate in $w$. Hence $w=(n-1)nw'(n-1)n$ where $w'$ is any 132-avoiding permutation over $[n-2]$. Thus, this case contributes $C_{n-2}$ representations.

\noindent
{\bf Subcase 1.3.}  $w_1=\epsilon$ and $w_3=1$. Similarly to the Subcases 1.1 and 1.2, we obtain that $w=nw'1$   where $w'$ is the increasing permutation $12\cdots n$, and this case contributes one representation.

\noindent
{\bf Subcase 1.4.}  $w_1=\epsilon$ and $w_3=\epsilon$. Here, $w=nw_2n$
    where $w_2$ is a 132-avoiding permutation over $[n-1]$. Thus, this case contributes $C_{n-1}$ representations.\\

\noindent
{\bf Case 2.} There is only one copy of $n$ in $w$. For $1\leq i< j\leq n-1$, suppose that there are exactly two copies of $i$ and $j$ in $w$ (by Theorem~\ref{degree2} there can be at most two copies of each letter). Since $K_n$ is a complete graph, we have that $n$ lies between $i_1$ and $i_2$ in $w$, and $n$ also lies between $j_1$ and $j_2$ in $w$, where recall that, e.g. $i_1$ and $i_2$ denote the first and the second occurrences of $i$, respectively, in the word. Then $i_1,n,j_2$ will form the pattern 132; contradiction. Using Theorem~\ref{degree2}, there are two subcases to consider.

\noindent
{\bf Subcase 2.1.} Every element in $A(w)$ has only one occurrence in $w$. Thus, $w$ is a 132-avoiding permutation over $\{1,2,\ldots,n\}$. Thus, this case contributes $C_{n}$ representations.

\noindent
{\bf Subcase 2.2.} There is only one letter $i$, $1\leq i \leq n-1$, in $A(w)$ that occurs twice in $w$.  Any letter  in $A(w)$ distinct from $i$ must lie between $i_1$ and $i_2$ in $w$. Since $w$ is 132-avoiding, we obtain that $w=i(i+1)\cdots nw'i$ where $w'$ is any 132-avoiding permutation over $\{1,2,\ldots,i-1\}$. Thus, this case contributes $\sum_{i=1}^{n-1}C_{i-1}=\sum_{i=0}^{n-2}C_{i}$ representations.
\qed

By Theorem~\ref{thm-K-n-132}, the initial values for the number of 132-representants for $K_n$, starting from $n=3$, are
$$12, 27, 72, 213, 670, 2190, 7349, 25146, 87364, 307310, 1092200, 3915866,\ldots.$$

\begin{exa}
\textnormal{
For $n=3$, we can see that all 12 132-representants for $K_3$, ordered as in the proof of Theorem~\ref{thm-K-n-132}, are 231231; 23123; 31231; 3123, 3213;~ 123, 231, 213, 312, 321; 1231, 2312.
}
\end{exa}

A direct corollary of Theorem~\ref{thm-K-n-132} is the following statement.

\begin{coro} For $n\geq 3$ and a $132$-representant $w$ for  $K_n$, the length of $w$ is either $n$, or $n+1$, or $n+2$, or $n+3$.
\end{coro}

\subsection{Non-132-representable graphs and 132-representation of small graphs}

Each non-word-representable graph is clearly non-132-representable. In this subsection we will show that the minimum (with respect to the number of vertices) non-word-representable graph, the wheel graph $W_5$ given in Figure~\ref{wheel5}, is actually a minimum non-132-representable graph. We do not know whether there exist other non-132-representable graphs on six vertices (no other non-word-representable graphs on six vertices exist). As for non-132-representable but word-representable graphs, an example of those is prisms $Pr_n$, where $n\geq 3$. The latter follows from Theorems~\ref{3-repr-not-2-repr} and~\ref{thm3-2}.

We note that the complement of a 132-representable graph is not necessarily a 132-representable graph. Indeed, for example, the 132-avoiding word 6645342312 defines a 132-representable graph, which is disjoint union of a cycle and the isolated vertex 6. However, the complement of this graph is the wheel graph $W_5$, which is not word-representable.

The following lemma allows us to restrict ourselves to considering graphs without isolated vertices when studying 132-representation.

\begin{lem}\label{no-isolat-added}
Let $G$ be a graph and $G'$ be a graph obtained from $G$ by adding an isolated vertex. Then $G$ is $132$-representable if and only if $G'$ is $132$-representable. \end{lem}

\proof
If $G'$ is 132-represented by $w$ then removing from $w$ the letter corresponding to  the isolated vertex we obtain a word 132-representing $G$.

Conversely, suppose that $G$ is 132-represented by $w$ and $n$ is larger than any letter in $w$. Then we label the isolated vertex by $n$ and note that the word $nnw$ 132-represents $G'$.
\qed

Lemma~\ref{no-isolat-added} cannot be generalized to adding to a graph a new connected 132-representable component instead of an isolated vertex. This follows from the fact established in \cite{M} that disjoint union of two complete graphs $K_4$ is non-132-representable, while $K_4$ is 132-representable. However, such a generalization can be done in a special case as recorded in the following simple, but useful lemma.

\begin{lem}\label{useful-lemma}
Let $G_1$, $G_2,\ldots,G_k$ be connected components of a graph $G$ that can be $132$-represented by  $2$-uniform words $w_1$, $w_2,\ldots,w_k$, respectively. Then $G$ is $132$-representable (by a $2$-uniform word).
\end{lem}

\proof
For $1\leq i\leq k$, let $a_i=|A(w_i)|$ denote the number of vertices in $G_i$, and let $\red^*(w_i)$ denote the word obtained from $\red(w_i)$ by replacing each element $j$, $1\leq j \leq a_i$, by $j+\sum_{m=1}^{i-1}a_m$. Then the 2-uniform word
$$w=\red^*(w_k)\red^*(w_{k-1})\cdots \red^*(w_1)$$
$132$-represents $G$.
\qed

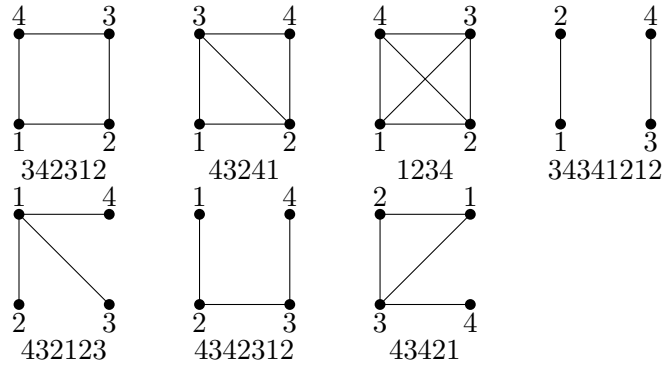
\begin{figure}[H]
    \centering
    \begin{minipage}{.75\textwidth}
        \centering
        \begin{tikzpicture}[scale=1.2]
       \draw (3,1)--(2,1)--(2,0);
       \draw (2,1)--(3,0);

       \fill(2,0) circle(0.06cm);\fill(2,1) circle(0.06cm);
       \fill(3,0) circle(0.06cm);\fill(3,1) circle(0.06cm);

       \node  at (2,-0.2) {  $2$};\node at (2,1.2) {  $1$};
       \node  at (3,-0.2) {  $3$};\node  at (3,1.2) {  $4$};

       \node  at (2.5,-0.5) {  $432123$};

       \draw (4,1)--(4,0)--(5,0)--(5,1);

       \fill(4,0) circle(0.06cm);\fill(4,1) circle(0.06cm);
       \fill(5,0) circle(0.06cm);\fill(5,1) circle(0.06cm);

       \node  at (4,-0.2) {  $2$};\node at (4,1.2) {  $1$};
       \node  at (5,-0.2) {  $3$};\node  at (5,1.2) {  $4$};

       \node  at (4.5,-0.5) {  $4342312$};

     \draw (7,1)--(6,1)--(6,0)--(7,0);
     \draw (6,0)--(7,1);

       \fill(6,0) circle(0.06cm);\fill(6,1) circle(0.06cm);
       \fill(7,0) circle(0.06cm);\fill(7,1) circle(0.06cm);

       \node  at (6,-0.2) {  $3$};\node at (6,1.2) {  $2$};
       \node  at (7,-0.2) {  $4$};\node  at (7,1.2) {  $1$};

       \node  at (6.5,-0.5) {  $43421$};

      \draw (2,2)--(3,2)--(3,3)--(2,3)--(2,2);

       \fill(2,2) circle(0.06cm);\fill(2,3) circle(0.06cm);
       \fill(3,2) circle(0.06cm);\fill(3,3) circle(0.06cm);

       \node  at (2,1.8) {  $1$};\node at (2,3.2) {  $4$};
       \node  at (3,1.8) {  $2$};\node  at (3,3.2) {  $3$};

       \node  at (2.5,1.5) {  $342312$};
       \draw (4,2)--(5,2)--(5,3)--(4,3)--(4,2);
       \draw (4,3)--(5,2);

       \fill(4,2) circle(0.06cm);\fill(4,3) circle(0.06cm);
       \fill(5,2) circle(0.06cm);\fill(5,3) circle(0.06cm);

       \node  at (4,1.8) {  $1$};\node at (4,3.2) {  $3$};
       \node  at (5,1.8) {  $2$};\node  at (5,3.2) {  $4$};

       \node  at (4.5,1.5) {  $43241$};
      \draw (6,2)--(7,2)--(7,3)--(6,3)--(6,2);
       \draw (6,3)--(7,2);
       \draw (6,2)--(7,3);

       \fill(6,2) circle(0.06cm);\fill(6,3) circle(0.06cm);
       \fill(7,2) circle(0.06cm);\fill(7,3) circle(0.06cm);

       \node  at (6,1.8) {  $1$};\node at (6,3.2) {  $4$};
       \node  at (7,1.8) {  $2$};\node  at (7,3.2) {  $3$};

       \node  at (6.5,1.5) {  $1234$};
       \draw (8,2)--(8,3);
       \draw (9,2)--(9,3);

       \fill(8,2) circle(0.06cm);\fill(8,3) circle(0.06cm);
       \fill(9,2) circle(0.06cm);\fill(9,3) circle(0.06cm);

       \node  at (8,1.8) {  $1$};\node at (8,3.2) {  $2$};
       \node  at (9,1.8) {  $3$};\node  at (9,3.2) {  $4$};

       \node  at (8.5,1.5) {  $34341212$};
   \end{tikzpicture}
    \end{minipage}
  \caption{$132$-representants for graphs on four vertices}
  \label{fig-4-132-repr}
\end{figure}

\begin{figure}[H]
    \centering
    \begin{minipage}{.95\textwidth}
        \centering
        \begin{tikzpicture}[scale=1]

       \draw (0,0)--(1,1);\draw (0,0)--(1,0);\draw (0,0)--(1.732,0.5);
       \draw (0,1)--(1,0);\draw (0,1)--(1,1);\draw (0,1)--(1.732,0.5);
       \draw (1,0)--(1.732,0.5);\draw (1,1)--(1.732,0.5);

       \fill(0,0) circle(0.06cm);\fill(1,0) circle(0.06cm);
       \fill(1,1) circle(0.06cm);\fill(0,1) circle(0.06cm);
       \fill(1.732,0.5) circle(0.06cm);

       \node  [left] at (0,-0.1) {  $2$};\node  [right] at (1,-0.1) {  $4$};
       \node  [left] at (0,1.1) {  $5$};\node   [right] at (1,1.1) {  $1$};
       \node [right] at (1.732,0.5) {  $3$};

       \node  at (0.866,-0.4) {  $45234125$};

       \draw (2.5,0)--(3.5,1);\draw (2.5,0)--(3.5,0);\draw (2.5,0)--(4.232,0.5);
       \draw (2.5,1)--(3.5,0);\draw (2.5,1)--(3.5,1);\draw (2.5,1)--(4.232,0.5);
       \draw (3.5,0)--(4.232,0.5);\draw (3.5,1)--(4.232,0.5);
       \draw (3.5,0)--(3.5,1);

       \fill(2.5,0) circle(0.06cm);\fill(3.5,0) circle(0.06cm);
       \fill(3.5,1) circle(0.06cm);\fill(2.5,1) circle(0.06cm);
       \fill(4.232,0.5) circle(0.06cm);

      \node [left] at (2.5,-0.1) {  $1$};\node  [right] at (3.5,-0.1) {  $3$};
       \node  [left] at (2.5,1.1) {  $5$};\node   [right] at (3.5,1.1) {  $4$};
       \node [right] at (4.232,0.5) {  $2$};

        \node  at (3.366,-0.4) {  $512341$};
       \draw (5,0)--(6,1);\draw (5,0)--(6,0);\draw (5,0)--(6.732,0.5);
       \draw (5,1)--(6,0);\draw (5,1)--(6,1);\draw (5,1)--(6.732,0.5);
       \draw (6,0)--(6.732,0.5);\draw (6,1)--(6.732,0.5);
       \draw (5,0)--(5,1);\draw (6,0)--(6,1);

       \fill(5,0) circle(0.06cm);\fill(6,0) circle(0.06cm);
       \fill(6,1) circle(0.06cm);\fill(5,1) circle(0.06cm);
       \fill(6.732,0.5) circle(0.06cm);

      \node [left] at (5,-0.1) {  $1$};\node  [right] at (6,-0.1) {  $2$};
       \node  [left] at (5,1.1) {  $3$};\node   [right] at (6,1.1) {  $4$};
       \node [right] at (6.732,0.5) {  $5$};

       \node  at (5.866,-0.4) {  $54321$};

        \draw (0,2)--(1,3);\draw (0,2)--(1,2);
       \draw (0,3)--(1,2);\draw (0,3)--(1,3);
       \draw (1,2)--(1.732,2.5);
       \draw (0,2)--(0,3);\draw (1,2)--(1,3);

       \fill(0,2) circle(0.06cm);\fill(1,2) circle(0.06cm);
       \fill(1,3) circle(0.06cm);\fill(0,3) circle(0.06cm);
       \fill(1.732,2.5) circle(0.06cm);

       \node  [left] at (0,1.9) {  $2$};\node  [right] at (1,1.9) {  $4$};
       \node  [left] at (0,3.1) {  $1$};\node   [right] at (1,3.1) {  $3$};
       \node [right] at (1.732,2.5) {  $5$};

       \node  at (0.866,1.6) {  $545321$};

        \draw (2.5,2)--(3.5,3);\draw (2.5,2)--(3.5,2);
       \draw (2.5,3)--(3.5,2);\draw (2.5,3)--(3.5,3);
       \draw (3.5,2)--(4.232,2.5);\draw (3.5,3)--(4.232,2.5);
       \draw (3.5,2)--(3.5,3);

       \fill(2.5,2) circle(0.06cm);\fill(3.5,2) circle(0.06cm);
       \fill(3.5,3) circle(0.06cm);\fill(2.5,3) circle(0.06cm);
       \fill(4.232,2.5) circle(0.06cm);

       \node [left] at (2.5,1.9) {  $4$};\node  [right] at (3.5,1.9) {  $3$};
       \node  [left] at (2.5,3.1) {  $1$};\node   [right] at (3.5,3.1) {  $2$};
       \node [right] at (4.232,2.5) {  $5$};

       \node  at (3.366,1.6) {  $5432451$};
         \draw (5,2)--(6,3);\draw (5,2)--(6,2);
       \draw (5,3)--(6,2);\draw (5,3)--(6,3);
       \draw (6,2)--(6.732,2.5);\draw (6,3)--(6.732,2.5);
       \draw (5,2)--(5,3);

       \fill(5,2) circle(0.06cm);\fill(6,2) circle(0.06cm);
       \fill(6,3) circle(0.06cm);\fill(5,3) circle(0.06cm);
       \fill(6.732,2.5) circle(0.06cm);

        \node [left] at (5,1.9) {  $3$};\node  [right] at (6,1.9) {  $5$};
       \node  [left] at (5,3.1) {  $4$};\node   [right] at (6,3.1) {  $1$};
       \node [right] at (6.732,2.5) {  $2$};

       \node  at (5.866,1.6) {  $5432512$};
       \draw (7.5,2)--(8.5,3);\draw (7.5,2)--(8.5,2);\draw (7.5,3)--(8.5,3);
       \draw (8.5,2)--(9.232,2.5);\draw (8.5,3)--(9.232,2.5);
       \draw (7.5,2)--(7.5,3);\draw (8.5,2)--(8.5,3);

       \fill(7.5,2) circle(0.06cm);\fill(8.5,2) circle(0.06cm);
       \fill(8.5,3) circle(0.06cm);\fill(7.5,3) circle(0.06cm);
       \fill(9.232,2.5) circle(0.06cm);

       \node [left] at (7.5,1.9) {  $4$};\node  [right] at (8.5,1.9) {  $5$};
       \node  [left] at (7.5,3.1) {  $1$};\node   [right] at (8.5,3.1) {  $3$};
       \node [right] at (9.232,2.5) {  $2$};

       \node  at (8.366,1.6) {  $5423521$};

       \draw (10,2)--(11,3);\draw (10,2)--(11,2);
       \draw (10,3)--(11,2);\draw (10,3)--(11,3);
       \draw (11,2)--(11.732,2.5);\draw (11,3)--(11.732,2.5);
       \draw (10,2)--(10,3);\draw (11,2)--(11,3);

       \fill(10,2) circle(0.06cm);\fill(11,2) circle(0.06cm);
       \fill(11,3) circle(0.06cm);\fill(10,3) circle(0.06cm);
       \fill(11.732,2.5) circle(0.06cm);

       \node [left] at (10,1.9) {  $2$};\node  [right] at (11,1.9) {  $3$};
       \node  [left] at (10,3.1) {  $1$};\node   [right] at (11,3.1) {  $4$};
       \node [right] at (11.732,2.5) {  $5$};

      \node  at (10.866,1.6) {  $543521$};
       \draw (0,4)--(1,4);
       \draw (0,5)--(1,4);\draw (0,5)--(1,5);
       \draw (1,5)--(1.732,4.5);
       \draw (0,4)--(0,5);\draw (1,4)--(1,5);

       \fill(0,4) circle(0.06cm);\fill(1,4) circle(0.06cm);
       \fill(1,5) circle(0.06cm);\fill(0,5) circle(0.06cm);
       \fill(1.732,4.5) circle(0.06cm);

       \node  [left] at (0,3.9) {  $1$};\node  [right] at (1,3.9) {  $2$};
       \node  [left] at (0,5.1) {  $3$};\node   [right] at (1,5.1) {  $4$};
       \node [right] at (1.732,4.5) {  $5$};

       \node  at (0.866,3.6) {  $5453241$};

        \draw (2.5,4)--(3.5,4);\draw (2.5,5)--(3.5,4);\draw (2.5,5)--(3.5,5);
        \draw (3.5,4)--(4.232,4.5);\draw (2.5,4)--(2.5,5);\draw (3.5,4)--(3.5,5);

       \fill(2.5,4) circle(0.06cm);\fill(3.5,5) circle(0.06cm);
       \fill(3.5,4) circle(0.06cm);\fill(2.5,5) circle(0.06cm);
       \fill(4.232,4.5) circle(0.06cm);

      \node [left] at (2.5,3.9) {  $1$};\node  [right] at (3.5,3.9) {  $3$};
       \node  [left] at (2.5,5.1) {  $2$};\node   [right] at (3.5,5.1) {  $5$};
       \node [right] at (4.232,4.5) {  $4$};

        \node  at (3.366,3.6) {  $5434251$};
       \draw (5,4)--(6,4);\draw (5,5)--(6,5);
       \draw (6,4)--(6.732,4.5);\draw (6,5)--(6.732,4.5);
       \draw (5,4)--(5,5);\draw (6,4)--(6,5);

       \fill(5,4) circle(0.06cm);\fill(6,4) circle(0.06cm);
       \fill(6,5) circle(0.06cm);\fill(5,5) circle(0.06cm);
       \fill(6.732,4.5) circle(0.06cm);

        \node [left] at (5,3.9) {  $2$};\node  [right] at (6,3.9) {  $1$};
       \node  [left] at (5,5.1) {  $3$};\node   [right] at (6,5.1) {  $4$};
       \node [right] at (6.732,4.5) {  $5$};

        \node  at (5.866,3.6) {  $5342312$};
       \draw (7.5,4)--(8.5,5);\draw (7.5,4)--(8.5,4);
       \draw (7.5,5)--(8.5,4);\draw (7.5,5)--(8.5,5);
       \draw (8.5,4)--(9.232,4.5);\draw (8.5,5)--(9.232,4.5);

       \fill(7.5,4) circle(0.06cm);\fill(8.5,4) circle(0.06cm);
       \fill(8.5,5) circle(0.06cm);\fill(7.5,5) circle(0.06cm);
       \fill(9.232,4.5) circle(0.06cm);

       \node [left] at (7.5,3.9) {  $5$};\node  [right] at (8.5,3.9) {  $1$};
       \node  [left] at (7.5,5.1) {  $4$};\node   [right] at (8.5,5.1) {  $3$};
       \node [right] at (9.232,4.5) {  $2$};

       \node  at (8.366,3.6) {  $34523125$};
        \draw (10,4)--(11,5);\draw (10,5)--(11,5);
       \draw (11,4)--(11.732,4.5);\draw (11,5)--(11.732,4.5);
       \draw (10,4)--(10,5);\draw (11,4)--(11,5);

       \fill(10,4) circle(0.06cm);\fill(11,4) circle(0.06cm);
       \fill(11,5) circle(0.06cm);\fill(10,5) circle(0.06cm);
       \fill(11.732,4.5) circle(0.06cm);

        \node [left] at (10,3.9) {  $1$};\node  [right] at (11,3.9) {  $3$};
       \node  [left] at (10,5.1) {  $2$};\node   [right] at (11,5.1) {  $4$};
       \node [right] at (11.732,4.5) {  $5$};

       \node  at (10.866,3.6) {  $5345312$};
     \draw (0,6)--(1,6);\draw (0,7)--(1,7);
      \draw (1,6)--(1.732,6.5);\draw (1,7)--(1.732,6.5);
      \draw (1,6)--(1,7);

       \fill(0,6) circle(0.06cm);\fill(1,6) circle(0.06cm);
       \fill(1,7) circle(0.06cm);\fill(0,7) circle(0.06cm);
       \fill(1.732,6.5) circle(0.06cm);

       \node  [left] at (0,5.9) {  $1$};\node  [right] at (1,5.9) {  $2$};
       \node  [left] at (0,7.1) {  $5$};\node   [right] at (1,7.1) {  $4$};
       \node [right] at (1.732,6.5) {  $3$};

       \node  at (0.866,5.6) {  $5453121$};

       \draw (2.5,6)--(3.5,7);\draw (2.5,7)--(3.5,7);
       \draw (3.5,6)--(4.232,6.5);\draw (3.5,7)--(4.232,6.5);
       \draw (3.5,6)--(3.5,7);

       \fill(2.5,6) circle(0.06cm);\fill(3.5,6) circle(0.06cm);
       \fill(3.5,7) circle(0.06cm);\fill(2.5,7) circle(0.06cm);
       \fill(4.232,6.5) circle(0.06cm);

       \node [left] at (2.5,5.9) {  $4$};\node  [right] at (3.5,5.9) {  $1$};
       \node  [left] at (2.5,7.1) {  $5$};\node   [right] at (3.5,7.1) {  $3$};
       \node [right] at (4.232,6.5) {  $2$};

      \node  at (3.366,5.6) {  $5434521$};
      \draw (5,7)--(6,7); \draw (6,6)--(6.732,6.5);\draw (6,7)--(6.732,6.5);
       \draw (5,6)--(5,7);\draw (6,6)--(6,7);

       \fill(5,6) circle(0.06cm);\fill(6,6) circle(0.06cm);
       \fill(6,7) circle(0.06cm);\fill(5,7) circle(0.06cm);
       \fill(6.732,6.5) circle(0.06cm);

        \node [left] at (5,5.9) {  $5$};\node  [right] at (6,5.9) {  $2$};
       \node  [left] at (5,7.1) {  $4$};\node   [right] at (6,7.1) {  $3$};
       \node [right] at (6.732,6.5) {  $1$};

    \node  at (5.866,5.6) {  $5453421$};
       \draw (7.5,6)--(8.5,6);\draw (7.5,7)--(8.5,7);
       \draw (8.5,6)--(9.232,6.5);
       \draw (7.5,6)--(7.5,7);\draw (8.5,6)--(8.5,7);

       \fill(7.5,6) circle(0.06cm);\fill(8.5,6) circle(0.06cm);
       \fill(8.5,7) circle(0.06cm);\fill(7.5,7) circle(0.06cm);
       \fill(9.232,6.5) circle(0.06cm);

       \node [left] at (7.5,5.9) {  $4$};\node  [right] at (8.5,5.9) {  $1$};
       \node  [left] at (7.5,7.1) {  $3$};\node   [right] at (8.5,7.1) {  $2$};
       \node [right] at (9.232,6.5) {  $5$};

       \node  at (8.366,5.6) {  $3423124515$};
       \draw (10,6)--(11,6) ;\draw (10,7) --(11,7);
       \draw (11,6) --(11.732,6.5);\draw (11,7)--(11.732,6.5);
       \draw (10,6)--(10,7) ;

       \fill(10,6) circle(0.06cm);\fill(11,6) circle(0.06cm);
       \fill(11,7) circle(0.06cm);\fill(10,7) circle(0.06cm);
       \fill(11.732,6.5) circle(0.06cm);

        \node [left] at (10,5.9) {  $1$};\node  [right] at (11,5.9) {  $2$};
       \node  [left] at (10,7.1) {  $5$};\node   [right] at (11,7.1) {  $4$};
       \node [right] at (11.732,6.5) {  $3$};

        \node  at (10.866,5.6) {  $45342312$};
      \draw (1,9)--(1.732,8.5);\draw (0,8)--(0,9);\draw (1,8)--(1,9);

       \fill(0,8) circle(0.06cm);\fill(1,8) circle(0.06cm);
       \fill(1,9) circle(0.06cm);\fill(0,9) circle(0.06cm);
       \fill(1.732,8.5) circle(0.06cm);

        \node  [left] at (0,7.9) {  $1$};\node  [right] at (1,7.9) {  $4$};
       \node  [left] at (0,9.1) {  $2$};\node   [right] at (1,9.1) {  $3$};
       \node [right] at (1.732,8.5) {  $5$};

        \node at (0.866,7.6) {  $543452121$};
       \draw (3.5,8)--(4.232,8.5);\draw (3.5,9)--(4.232,8.5);
       \draw (2.5,8)--(2.5,9);\draw (3.5,8)--(3.5,9);

       \fill(2.5,8) circle(0.06cm);\fill(3.5,8) circle(0.06cm);
       \fill(3.5,9) circle(0.06cm);\fill(2.5,9) circle(0.06cm);
       \fill(4.232,8.5) circle(0.06cm);

      \node [left] at (2.5,7.9) {  $4$};\node  [right] at (3.5,7.9) {  $3$};
       \node  [left] at (2.5,9.1) {  $5$};\node   [right] at (3.5,9.1) {  $2$};
       \node [right] at (4.232,8.5) {  $1$};

        \node at (3.366,7.6) {  $5454321$};

       \draw (5,9)--(6,9);
       \draw (6,8)--(6.732,8.5);\draw (6,9)--(6.732,8.5);
       \draw (5,8)--(5,9);

       \fill(5,8) circle(0.06cm);\fill(6,8) circle(0.06cm);
       \fill(6,9) circle(0.06cm);\fill(5,9) circle(0.06cm);
       \fill(6.732,8.5) circle(0.06cm);

        \node [left] at (5,7.9) {  $1$};\node  [right] at (6,7.9) {  $5$};
       \node  [left] at (5,9.1) {  $2$};\node   [right] at (6,9.1) {  $3$};
       \node [right] at (6.732,8.5) {  $4$};

     \node at (5.866,7.6) {  $545342312$};
       \draw (7.5,9)--(8.5,9);\draw (8.5,9)--(9.232,8.5);
       \draw (7.5,8)--(7.5,9);\draw (8.5,8)--(8.5,9);

       \fill(7.5,8) circle(0.06cm);\fill(8.5,8) circle(0.06cm);
       \fill(8.5,9) circle(0.06cm);\fill(7.5,9) circle(0.06cm);
       \fill(9.232,8.5) circle(0.06cm);

       \node [left] at (7.5,7.9) {  $1$};\node  [right] at (8.5,7.9) {  $4$};
       \node  [left] at (7.5,9.1) {  $2$};\node   [right] at (8.5,9.1) {  $3$};
       \node [right] at (9.232,8.5) {  $5$};

        \node at (8.366,7.6) {  $543452312$};
       \draw (10,8) --(11.732,8.5);\draw (10,9)--(11.732,8.5);
       \draw (11,8)--(11.732,8.5);\draw (11,9)--(11.732,8.5);

       \fill(10,8) circle(0.06cm);\fill(11,8) circle(0.06cm);
       \fill(11,9) circle(0.06cm);\fill(10,9) circle(0.06cm);
       \fill(11.732,8.5) circle(0.06cm);

        \node [left] at (10,7.9) {  $4$};\node  [right] at (11,7.9) {  $5$};
       \node  [left] at (10,9.1) {  $3$};\node   [right] at (11,9.1) {  $2$};
       \node [right] at (11.732,8.5) {  $1$};

        \node at (10.866,7.6) {  $543212345$};


        \end{tikzpicture}
    \end{minipage}
  \caption{$132$-representants for graphs on five vertices (with no isolated vertices)}
  \label{graphs-5}
\end{figure}
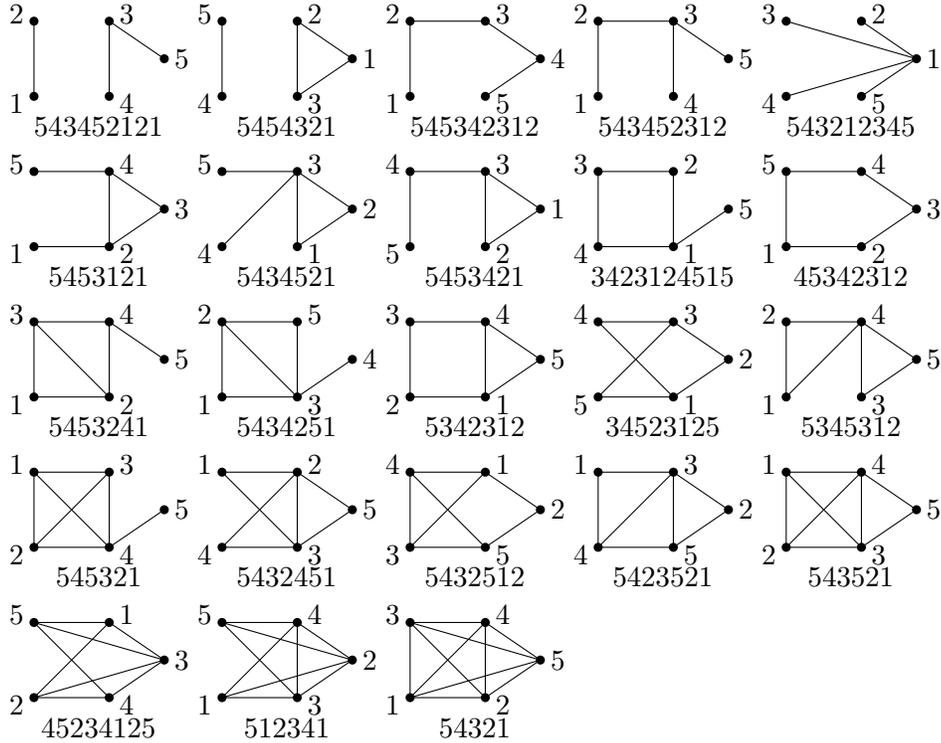

By Lemma~\ref{no-isolat-added}, we exclude isolated vertices from our considerations in the rest of this subsection.  Moreover, graphs on up to three vertices are either trees or the cycle graph $C_3$, and thus they are 132-representable. Further, there are seven graphs on four vertices which can be 132-represented as shown in Figure~\ref{fig-4-132-repr}. Finally, there are $23$  graphs on five vertices that have no isolated vertices, and these graphs can be 132-represented as in Figure~\ref{graphs-5}. Note that Lemma~\ref{useful-lemma} was used (in a straightforward way) to 132-represent graphs in Figures~\ref{fig-4-132-repr} and~\ref{graphs-5} that have two connected components.

\section{Concluding remarks}\label{open-sec}

This paper just scratches the surface of a big research direction dealing with representing graphs by pattern-avoiding words. Our studies were extended to 123-representation of graphs in \cite{M}, where more results on 132-representable graphs were obtained as well. Further steps may be in considering longer patterns and/or patterns of other types (e.g. those described in~\cite{HM,Kit}) while defining words to be used to represent graphs, and asking the question on which classes of graphs can be represented in this way. Simultaneous avoidance of patterns, like avoiding the patterns 132 and 231 at the same time, can be considered as well.

To conclude, we state the following question, solving which by exhaustive search would involve finding appropriate labelling of graphs and then considering all words over six letter alphabet that have at most two occurrences of each letter.

\noindent
\textbf{Question:} Is the wheel graph $W_5$ the only non-132-representable graph on six vertices?

\section*{Acknowledgments}

The work of the first and the third authors was supported by the 973 Project, the PCSIRT Project of the Ministry of Education and the National Science Foundation of China. The  second author is grateful to Bill Chen and Arthur Yang for their hospitality during the author's stay at the Center for Combinatorics at Nankai University in November 2015. All the authors are also grateful to the  Center for Applied Mathematics at Tianjin University for its generous support.


\begin{thebibliography}{13}

%
%
%
%
%
%
%
\bibitem{Halldorsson15} M. Halld\'{o}rsson, S. Kitaev, and A. Pyatkin. Semi-transitive orientations and word-representable graphs, {\em Discr. Appl. Math.}  {\bf 201} (2016) 164--171.

\bibitem{HM} S. Heubach and T. Mansour. {\em Combinatorics of compositions
   and words}, Chapman \& Hall/CRC an imprint of Taylor \& Francis LLC,
   Discrete Mathematics and Its Applications Series, 2009.

\bibitem{JKPR} M. Jones, S. Kitaev, A. Pyatkin and J. Remmel. Representing graphs via pattern avoiding words, {\em Elect. J. Combin.} {\bf 22(2)} (2015), \#P2.53, 20 pp.

\bibitem{Kit} S. Kitaev. {\em Patterns in Permutations and Words}, Springer, 2011.

\bibitem{Kit1} S. Kitaev. On graphs with representation number 3, {\em J. Autom. Lang. Comb.} {\bf 18} (2013) 2, 97--112.

\bibitem{KL} S. Kitaev and V. Lozin. {\em Words and Graphs}, Springer, 2015.

\bibitem{KP} S. Kitaev and A. Pyatkin. On representable graphs {\em J. Autom. Lang. Comb.} {\bf 13} (2008) 1, 45--54.
%
%
%

\bibitem{M} Y. Mandelshtam. On graphs representable by pattern-avoiding words, 	arXiv:1608.07614.

\end{thebibliography}
\end{document}